\documentclass[12pt]{article}
\usepackage{amsmath}
\usepackage{amssymb}
\usepackage{latexsym}
\makeatletter
   \renewcommand{\section}{\@startsection {section}{1}{\z@}%
   {-3.5ex \@plus -1ex \@minus -.2ex}%
   {2.3ex \@plus.2ex}%
   {\normalfont\large\bfseries}}
\makeatother

\newcommand{\hovedfont}{\normalfont\bfseries}
\usepackage{theorem}
\theorembodyfont{\normalfont\itshape}
\theoremheaderfont{\hovedfont}
   \theoremstyle{change}
\newtheorem{lemma}{Lemma.}[section]
\newtheorem{satz}[lemma]{Theorem.}

   \theorembodyfont{\normalfont}
\newtheorem{eks}[lemma]{Example.}
\newtheorem{obs}[lemma]{Remark.}
\newtheorem{taller}[lemma]{$\!\!$}
   \newenvironment{blanko}[1]%
{\begin{taller}{\hovedfont #1}\normalfont}%
{\end{taller}}
{\begin{list}{\textbf{Definition. }}%
{\setlength{\labelsep}{0mm}\setlength{\leftmargin}{0mm}%
\setlength{\labelwidth}{0mm}\setlength{\listparindent}{\parindent}%
\setlength{\parsep}{\parskip}\setlength{\partopsep}{0mm}}%
\item}{\end{list}}
   \newenvironment{dem}%
{\begin{list}{\em Proof. }%
{\setlength{\labelsep}{0mm}\setlength{\leftmargin}{0mm}%
\setlength{\labelwidth}{0mm}\setlength{\listparindent}{\parindent}%
\setlength{\parsep}{\parskip}\setlength{\partopsep}{0mm}}%
\item}{\qed\end{list}}
   \newenvironment{dem*}[1]%
{\begin{list}{\em #1 }%
{\setlength{\labelsep}{0mm}\setlength{\leftmargin}{0mm}%
\setlength{\labelwidth}{0mm}\setlength{\listparindent}{\parindent}%
\setlength{\parsep}{\parskip}\setlength{\partopsep}{0mm}}%
\item}{\qed\end{list}}
   \newenvironment{blanko*}[1]%
{\begin{list}{\bf {#1} }%
{\setlength{\labelsep}{0mm}\setlength{\leftmargin}{0mm}%
\setlength{\labelwidth}{0mm}\setlength{\listparindent}{\parindent}%
\setlength{\parsep}{\parskip}\setlength{\partopsep}{0mm}}%
\item}{\end{list}}
\renewcommand{\epsilon}{\varepsilon}
\makeatletter
\renewcommand{\ldots}{\relax\ifmmode\ldotp\ldotp\ldotp%
\else$\m@th\ldotp\ldotp\ldotp\ $\fi}
\makeatother
\newcommand{\brac}[1]{\,\langle \, {#1} \, \rangle\,}

\newcommand{\Brac}[1]{\,\langle \ {#1} \ \rangle\,}

\newcommand{\cp}{ {\scriptstyle \;\cup\;} }

\newcommand{\smallbinom}[2]{{\textstyle \binom{#1}{#2}}}
\newcommand{\smallprod}[2]{\overset{#2}{\underset{#1}{\textstyle{\prod}}}} 
\newcommand{\smallsum}[2]{\overset{#2}{\underset{#1}{\textstyle{\sum}}}} 

\newcommand{\bigsum}[2]{\overset{#2}{\underset{#1}{\text{{\Large\(\sum\)}}}}}

\newcommand{\bigfatgreek}[1]{\boldsymbol{#1}}

\newcommand{\psiclass}{\bigfatgreek\psi} 
\newcommand{\mpsiclass}{\ov\psiclass{}}

\newcommand{\deltaclass}{\bigfatgreek\delta} 
\newcommand{\boundary}{\bigfatgreek\xi}

\providecommand{\qed}{\hspace*{\fill}\nolinebreak[1]\hspace*{\fill} $\Box$}

\renewcommand{\d}{\partial}

\newcommand{\pil}       {\rightarrow}
\newcommand{\langpil}      {\longrightarrow}

\newcommand{\kan}          {\mbox{\Large \(\omega\)}}
\newcommand{\OO}{\ensuremath{{\mathcal O}}}

\newcommand{\upperstar}{^{\raisebox{-0.25ex}[0ex][0ex]{\(\ast\)}}}
\newcommand{\lowerstar}{_{\raisebox{-0.33ex}[-0.5ex][0ex]{\(\ast\)}}}
\newcommand{\df}{\: {\raisebox{0.255ex}{\normalfont\scriptsize :\!\!}}=}
\newcommand{\tensor} {\otimes}

\newcommand{\brok}[2]{{\textstyle{\frac{#1}{#2}}}}

\newcommand{\ov}{\overline}

\newcommand{\fat}[1]{\mathbf{{#1}}}

\renewcommand{\P} {\mathbb{P}}

\newcommand{\Z}      {\mathbb{Z}}
\newcommand{\Q}      {\mathbb{Q}}

\newcommand{\C}      {\mathbb{C}}

\newcommand{\mtau}{\ov\tau}
\newcommand{\mbtau}[2]{\ov{\boldsymbol\tau}_{\!#1}^{#2}}
\newcommand{\vdim}{\operatorname{vdim}}

\newcommand{\coh}[1]{\text{\sf\bfseries{{#1}}}}
\newcommand{\basis}{T}

\newcommand{\otrace}[1]{\int_{\coh 0}{#1}}
\newcommand{\ytrace}[1]{\int_{\coh y}{#1}}

\begin{document}

\vspace*{0.5in}
\begin{center}
   {\LARGE Tangency quantum cohomology \par}%
   \vskip 1.5em
   {\large \lineskip .5em
      Joachim Kock\footnote{%
      Supported by the National Science Research Council of Denmark.
      
      Date: 10th of June 2000.
      
      Current address: Dept.\ of Mathematics, \ Royal Institute of Technology, \
      100 44 Stockholm, Sweden.
         
      E-mail: {\tt kock@math.kth.se}}
   \par}%
   \vskip 1em
   {\normalsize Universidade Federal de Pernambuco\\
   Recife, Brazil}%
\end{center}%
\par
\vskip 1.5em

\begin{abstract}
   Let $X$ be a smooth projective variety.  Using modified psi classes on the
   stack of genus zero stable maps to $X$, a new associative quantum product is
   constructed on the cohomology space of $X$.  When $X$ is a homogeneous
   variety, this structure encodes the characteristic numbers of rational
   curves in $X$, and specialises to the usual quantum product upon resetting
   the parameters corresponding to the modified psi classes.  For $X=\P^2$, the
   product is equivalent to that of the contact cohomology of
   Ernstr\"om-Kennedy.
\end{abstract}

\section*{Introduction}

Let $X$ be a smooth projective variety over the complex numbers.  The
Gromov-Witten invariants of $X$ are constructed by pulling back cohomology
classes to the stack of stable maps to $X$ and integrating over the virtual
fundamental class (cf.~Behrend-Manin~\cite{Behrend-Manin},
Behrend~\cite{Behrend}, Li-Tian~\cite{Li-Tian}).  Let $\Phi$ denote the
generating function for the genus zero Gromov-Witten invariants.  The quantum
product is defined on the cohomology space of $X$ by taking the third
derivatives of $\Phi$ as structure constants,
$$
\basis_i * \basis_j \df \sum_{e,f} \Phi_{ije} \; g^{ef} \; \basis_f,
$$
cf.~Kontsevich-Manin~\cite{KM:9402}.  The associativity of this product is
equivalent to the fact that $\Phi$ satisfies the WDVV equations
\begin{eqnarray*}
	\sum_{e,f} \Phi_{ije} \; g^{ef} \; \Phi_{fkl} 
	&=& \pm
\sum_{e,f} \Phi_{jke} \; g^{ef} \; \Phi_{fil}.
\end{eqnarray*}
A feature of this structure is that when $X$ is a homogeneous variety, the
Gromov-Witten invariants are solutions to enumerative problems of rational
curves in $X$ subject to incidence conditions, and then the WDVV equations
provide relations among the numbers.

It is straightforward to generalise the construction by including psi classes in
the integrals, to get a larger quantum product: instead of using just the
primary Gromov-Witten potential as above, one uses the whole gravitational
Gromov-Witten potential.  This product is associative too, due to the
(generalised) WDVV equations (Witten~\cite{Witten}).  This larger product has
not brought about as much interest as the first one, for two reasons: the
gravitational descendants are not directly interpretable in enumerative
geometry, and the equations corresponding to associativity are ``perpendicular
to the descent'' in the sense that each equation relates only invariants with
the same psi classes as factors.

\bigskip

The present work introduces a new extension of the quantum product, the {\em
tangency quantum product}, and proves its associativity.  When $X$ is a
homogeneous variety, this structure encodes all the characteristic numbers of
rational curves in $X$.

The construction is based on a different potential $\Gamma$, called the {\em
tangency quantum potential}, which incorporates one level of {\em modified psi
classes}.  Modified psi classes are boundary corrections of the tautological psi
classes motivated from enumerative geometry: tangency conditions are easily
expressible in terms of modified psi classes.  They were first discovered and
studied by R.~Pandharipande (unpublished \cite{Pand:mail}) who used them to
establish a topological recursion for the characteristic numbers of rational
curves in $\P^2$, and later, independently, by myself~\cite{_9902}.  (This work
is subsumed in \mbox{Graber-Kock-Pandharipande~\cite{Grab-Kock-Pand}}).

In contrast to the tautological psi classes, the modified ones restrict to the
boundary in a non-trivial way, giving rise to certain {\em diagonal classes}. 
Therefore the splitting lemma (\ref{split2}), which is the engine in the proof
of the associativity, is complicated by correction terms.  This mess is handled
by the introduction of a new ``metric'' $\gamma_{ef}$ on the cohomology space of
$X$.  It is a deformation of the Poincar\'e metric $g_{ef}$ in much the same way
as the tangency quantum potential is a deformation of the Gromov-Witten
potential. Amazingly it encodes all the combinatorics of the diagonal
corrections in the splitting lemma, leading to certain ``deformed'' WDVV
equations.  These in turn are equivalent to associativity of the {\em tangency
quantum product}, defined as
$$
\basis_i * \basis_j \df \basis_i \cp \basis_j  + 
\sum_{e,f} \Gamma_{ije} \, \gamma^{ef} \, \basis_f.
$$
(This deformed metric also clarifies the topological recursion, as explained in
\cite{Grab-Kock-Pand}.)

The last section is devoted to a technical detail: since modified psi classes
are not defined in degree zero, the same goes for the tangency potential.  The
missing degree zero part is now constructed separately to account for the
classical part of the new product.  This construction shows that the new product
is in fact ``integrable'' in the sense that its structure constants (with
respect to $\gamma$) are the third derivatives of a single function.  With this
last piece of data, the results of the paper can be summarised in saying that
the new structure is that of a formal Frobenius (super)manifold (over the power
series ring in the new variables).

Acknowledgements: this work is a part of my PhD thesis~\cite{_thesis}, and I am
indebted to my advisor Israel Vainsencher.  This particular part was carried out
while I was visiting the California Institute of Technology; I am thankful to
Caltech for its hospitality, and in particular to Rahul Pandharipande for
encouraging my research, and for many stimulating conversations.  I
have also benefited from conversations with Marzia Polito and correspondence 
with Gary Kennedy.

\section{First enumerative descendants}

\begin{blanko}{Set-up.}  
	Throughout we work over the field of complex numbers.  Let $X$ denote a
	smooth projective variety, and let $\basis_0,\ldots,\basis_r$ denote the
	elements of a homogeneous basis of the cohomology space $H\upperstar (X,\Q)$. 
	As in Manin's book~\cite{Manin} we consider $H\upperstar (X,\Q)$ as a linear
	supermanifold.
   
	Let $\ov M_{0,S}(X,\beta)$ denote the moduli stack of Kontsevich stable maps
	of genus zero whose direct image in $X$ is of class $\beta\in H_2^+(X,\Z)$,
	and whose marking set is $S=\{p_1,\ldots,p_n\}$.
	For each mark $p_i$, let $\nu_i : \ov M_{0,S}(X,\beta) \pil X$ denote the
	evaluation morphism that sends the class of a map $\mu$ to $\mu(p_i)$;
	pull-backs along evaluation morphisms of cohomology classes of $X$ are
	referred to as {\em evaluation classes}.  The reader is referred to
	Fulton-Pandharipande~\cite{FP-notes} for definitions and basic properties of
	stable maps, Gromov-Witten invariants and quantum cohomology.
\end{blanko}

\begin{blanko}{Modified psi classes (cf.~\cite{Grab-Kock-Pand}).}
	Let $\pi_0 : \ov M_{0,S\cup\{p_0\}}(X,\beta) \pil \ov M_{0,S}(X,\beta)$
	denote the forgetful morphism that forgets the extra mark $p_0$.  Together
	with the sections $\sigma_i$ corresponding to each of the marks in $S$, it
	constitutes the universal family.  For each mark there is a psi class
	defined as $\psiclass_i \df c_1(\sigma_i\upperstar \kan_{\pi_0})$, where
	$\kan_{\pi_0}$ is the relative dualising sheaf of $\pi_0$.  When $\beta\neq
	0$, the {\em modified psi class} is $\mpsiclass_i \df\hat\pi_i\upperstar
	\psiclass_i$, where $\hat\pi_i : \ov M_{0,S}(X,\beta) \pil \ov
	M_{0,\{p_i\}}(X,\beta)$ is the forgetful morphism that forgets all marks but
	$p_i$.

	The importance of modified psi classes comes from enumerative geometry: for
	example, if $X$ is a homogeneous variety and $Z\subset X$ is a very ample
	hypersurface of degree $\coh z\in H^2(X,\Q)$, then the cycle of maps tangent
	to $Z$ at $p_i$ is of class $\nu_i\upperstar (\coh z) \cp ( \mpsiclass_i +
	\nu_i\upperstar (\coh z))$.  It then follows from a transversality argument
	that the characteristic numbers of rational curves are top products of
	evaluation classes and modified psi classes.  For this reason, these top
	products (integrated over the virtual fundamental class of $\ov
	M_{0,S}(X,\beta)$),
   $$
   \Brac{\mtau_{m_1}(\coh z_1)\cdots \mtau_{m_n}(\coh z_n)}\!_{\beta} \df
\int \ov\psiclass{}_1^{m_1}\cp\nu_1\upperstar (\coh z_1) \cp\cdots \cp
\ov\psiclass{}_n^{m_n}\cp\nu_n\upperstar (\coh z_n)
\ \cap [\ov M_{0,S}(X,\beta) ]^{\text{virt}}
   $$
	are called {\em enumerative descendants}.

The enumerative descendants that appear in the characteristic number problem are
only the {\em first} enumerative descendants, i.e.\ those made up of factors of
type $\mtau_0(\coh x)$ and $\mtau_1(\coh y)$, for $\coh x, \coh y \in
H\upperstar (X,\Q)$.  For simplicity, we will work in the fixed basis
$\basis_0,\ldots,\basis_r$ and employ the following notation for the
corresponding integrals:
\begin{equation}\label{tauab}
\brac{\mbtau{0}{\fat a} \ \mbtau{1}{\fat b}}\!_{\beta} \df 
\Brac{ \smallprod{k=0}{r} (\mtau_0(\basis_k))^{a_k}
\smallprod{\ell=0}{r} (\mtau_1(\basis_\ell))^{b_\ell}}\!_{\beta},
\end{equation}
where $\fat{a}=(a_0,\ldots,a_r)$ and $\fat{b}=(b_0,\ldots,b_r)$ are vectors of
non-negative integers.  The integral is zero unless $\sum_k
\brok{1}{2}\deg(\basis_k) (a_k+b_k+1) = \vdim_{\C} \ov M_{0,S}(X,\beta)$, where
$S$ is of cardinality $n= \sum (a_k + b_k)$.
\end{blanko}

\begin{blanko}{Diagonal classes (cf.~\cite{Grab-Kock-Pand}).}\label{diagonal}
	The {\em diagonal class} $\deltaclass_{ij} \in H^2(\ov M_{0,S}(X,\beta), \Q)$
	is defined as the sum of all boundary divisors having $p_i$ and $p_j$
	together on a contracting twig.  It can also be described as the pull-back
	from $\ov M_{0,\{p_i,p_j\}}(X,\beta)$ of the divisor $D_{ij}$ having both
	marks on a contracting twig; from this description it is clear that diagonal
	classes are invariant under pull-back along forgetful morphisms.

   Let $\pi_0 : \ov M_{0,S\cup\{p_0\}}(X,\beta) \pil \ov M_{0,S}(X,\beta)$
   denote the morphism that forgets $p_0$.  Then
   \begin{equation}\label{push-delta}
      \pi_0{}\lowerstar \deltaclass_{0i} = 1,
      \end{equation}
   (the virtual fundamental class of $\ov M_{0,S}(X,\beta)$).

   The diagonal classes owe their name to  the following properties.
   \begin{eqnarray}
      \deltaclass_{ij}\deltaclass_{ik} & = & 
      \deltaclass_{ij}\deltaclass_{jk}\notag\\
      -\deltaclass_{ij}^2 \; = \;\deltaclass_{ij}\mpsiclass_i  & = & 
      \deltaclass_{ij}\mpsiclass_j\notag\\
      \deltaclass_{ij}\nu_i\upperstar (\coh z) 
      & = & \deltaclass_{ij}\nu_j\upperstar (\coh z),\label{swap}
   \end{eqnarray}%
   for $\coh z\in H\upperstar (X,\Q)$.
\label{D_ij}
	The swapping property also holds when $\deltaclass_{ij}$ is replaced with
	$D_{ij}$.
\end{blanko}

\begin{blanko}{Key formula.}
   Let $D = D(S', \beta' \mid S'', \beta'')$ denote the boundary divisor which
   is the image of the morphism
$$
\rho_D : \ov M_{0,S' \cup \{x'\}}(X,\beta') 
\times_X \ov M_{0,S'' \cup \{x''\}}(X,\beta'')
\langpil \ov M_{0,S}(X,\beta)
$$
	consisting in gluing together the two maps at $x'$ and $x''$.  We will assume
	that $S'$ and $S''$ are non-empty and that $\beta'$ and $\beta''$ are
	non-zero.  Then $\rho_D$ is birational onto $D$.  For short, let $\ov M'$ and
	$\ov M''$ denote the factors of the fibred product above, and let
   $$
   \jmath_D : \ov M' \times_X \ov M'' \langpil \ov M' \times \ov M''
   $$
   denote the inclusion in the cartesian product.

   Now there is the following formula for restricting a modified psi class to the 
   boundary (assuming $p_i \in S'$).
   \begin{equation} \label{restr-mpsi}
   \rho_D\upperstar \ov \psiclass_i = 
\jmath_D\upperstar\big( \mpsiclass_i+\deltaclass_{ix'}\big).
   \end{equation} 
\end{blanko}

\begin{lemma}\label{split2}\textsc{Splitting lemma for first enumerative
	descendants.} Let $D= D(S',\beta' | S'',\beta'')$ be a boundary divisor with
	$\beta'$ and $\beta''$ non-zero.  Let $\coh x$ and $\coh y$ be generic even
	elements of $H\upperstar (X,\Q)$, and let there be given a product
	$\mtau_{0}(\coh x)^a\,\mtau_1(\coh y)^b$ such that the marks
	corresponding to $a' \leq a$ and $b' \leq b$ belong to $S'$, and the
	remaining marks, corresponding to $a'' \leq a$ and $b'' \leq b$, belong to
	$S''$.  Then \setlength{\multlinegap}{0pt}
   \begin{multline*}
      \brac{ D \cdot
   \mtau_{0}(\coh x)^a\,\mtau_1(\coh y)^b}\!_\beta =
   \sum_{p,q} \sum_{{s}',{s}''}
   \smallbinom{{b}'}{{s}'}\smallbinom{{b}''}{{s}''}
   \brac{\mtau_{0}(\coh x)^{a'}\,\mtau_1(\coh y)^{b'-s'} \,
   \mtau_0(\coh y^{{s}'}\!\cp \basis_p)}\!_{\beta'}
   \\ g^{pq} \  
   \brac{\mtau_0(\basis_q\cp\coh y^{{s}''})
   \,\mtau_{0}(\coh x)^{a''}\,\mtau_1(\coh y)^{b''-s''}}\!_{\beta''}.
   \end{multline*}
	{\rm The outer sum is over the splitting indices $p$ and $q$ running from $0$
	to $r$, and the inner sum is over all ${s}'$ and ${s}''$: the symbol
	$\mtau_{1}(\coh y)^{ b'- s'}$ makes sense only for ${s}'\leq{b}'$, but since
	there is a binomial coefficient $\binom{{b}'}{{s}'}$ in front of it, which is
	zero unless ${s}'\leq{b}'$ we will allow any ${s}'$.  Ditto for ${s}''$.  }
   \end{lemma}
\begin{dem}
   There are four ingredients in the proof of this formula: the first is of
   course the splitting axiom of Gromov-Witten theory
   (cf.~Kontsevich-Manin~\cite{KM:9402}), which accounts for the overall shape
   of the formula.  Second, formula~(\ref{restr-mpsi}) tells how each factor
   $\mtau_1(\coh y)$ (say realised at mark $p_i$ as $\mpsiclass_i \cp
   \nu_i\upperstar (\coh y)$) restricts to give $(\mpsiclass_i +
   \deltaclass_{ix})\cp\nu_i\upperstar (\coh y)$ on the twig containing $p_i$ ---
   here $x$ denotes the gluing mark of that twig.  Now (on each twig) expand the
   product of all these restrictions into a sum over $s$ (which accounts
   for the binomial coefficients).  Third, apply formula~(\ref{swap}) to write $
   \deltaclass_{ix}\cp \nu_i\upperstar (\coh y) = \deltaclass_{ix}\cp
   \nu_{x}\upperstar (\coh y)$, with the effect of accumulating evaluation classes
   on the gluing mark.  Finally, by equation~(\ref{push-delta}) we can push down
   each term involving a diagonal class $\deltaclass_{ix}$ along the forgetful
   morphism forgetting $p_i$; the effect is simply deleting $\deltaclass_{ix}$.
\end{dem}

\newcommand{\ring}{\Lambda}

\begin{blanko}{The tangency quantum potential }\label{Gamma}
	$\Gamma$ shall now be defined as the generating function for the first
	enumerative descendants.   
	Let $\coh x= \sum x_i \basis_i$ and $\coh y = \sum y_i \basis_i$ be
	generic even elements of $H\upperstar (X,\Q)$, and set
\begin{eqnarray}\label{Gammadef}
\Gamma(\coh{x},\coh{y}) 
&\df& \sum_{\beta>0} q^\beta 
\brac{ \exp( \mtau_0(\coh x) + \mtau_1(\coh y))}\!_\beta \\
&=& \sum_{\beta>0} q^\beta \sum_{a,b}
\brac{\frac{\mtau_0(\coh x)^a}{a!} \frac{\mtau_1(\coh y)^b}{b!}}\!_\beta \notag .
\end{eqnarray}
The inner sum is over {\em all} non-negative integers $a$ and $b$ --- this is
meaningful since there is no $\beta=0$ term in the outer sum.  (The degree-zero
case is considered and included below, in \ref{beta0} and \ref{interpret0}.) 
The coefficients $q^\beta$ are necessary only to ensure formal convergence; they
belong to the Novikov ring $\ring$ (as defined in Getzler~\cite{Getzler:9612}),
which will be our coefficient ring when treating the tangency quantum potential.

For convenience we identify $\coh x$ and $\coh y$ with their coordinate vectors
$\fat x = (x_0,\ldots,x_r)$ and $\fat y = (y_0,\ldots,y_r)$ with respect to the
basis $\basis_0,\ldots,\basis_r$.  The potential thus belongs to the
power series ring $\ring[[\fat x, \fat
y]]=\ring[[x_0,\ldots,x_r,y_0,\ldots,y_r]]$, and expands to
\begin{equation}\label{ybxa}
\Gamma(\fat{x},\fat{y})  = \sum_{\beta>0}q^\beta\sum_{\fat{a},\fat{b}} 
\frac{\fat{y}^{\fat{b}}}{\fat{b}!}\frac{\fat{x}^{\fat{a}}}{\fat{a}!}
\brac{\mbtau{0}{\fat a} \ \mbtau{1}{\fat b}}\!_{\beta}.
\end{equation}
Here
the inner sum is over all pairs of vectors
$\fat{a}=(a_0,\ldots,a_r)$ and $\fat{b}=(b_0,\ldots,b_r)$ of non-negative
integers, and we employ multi-index notation,
e.g. $\fat a!= a_0!\cdots a_r!$.  For the formal variables $\fat x$ and $\fat y$,
the multi-index notation is reverse, in order to preserve the signs arising from 
odd variables, e.g.\
$\fat x^{\fat a} = x_r^{a_r}\cdots x_0^{a_0}$

The variables $\fat{x}$ are the usual formal variables from quantum cohomology,
so when $\fat{y}$ is set to zero, $\Gamma$ reduces to the usual (quantum part of
the) genus-zero Gromov-Witten potential.

The viewpoint of (\ref{ybxa}) is advantageous for the sake of extracting the
invariants, as well as for checking the validity of certain formal operations on
the potential.  This task is safely left to the reader and henceforth only the
more compact notation of (\ref{Gammadef}) is used.
\end{blanko}

\section{Deformation of the Poincar\'e metric}

While the usual quantum potential is based on the Poincar\'e metric constants
$g_{ij} = \otrace{\basis_i\cp\basis_j}$, the tangency quantum potential relates
more naturally to a deformation of them, a certain ``metric'' with
values in $\Q[[\fat{y}]]=\Q[[y_0,\ldots,y_r]]$.

\begin{blanko}{The classical product in the Poincar\'e metric.}
	The two important structures on $H=H\upperstar (X,\Q)$ are the intersection
	product $\cp$, and the trace map $\otrace{} : H \pil \Q$ which is just
	integration over the fundamental class of $X$.  Set
	$g_{ij}=\otrace{\basis_i\cp\basis_j}$ and $g_{ijk} = \otrace{\basis_i \cp
	\basis_j \cp \basis_k}$.  Let $(g^{ij})$ be the inverse matrix to $(g_{ij})$. 
	It is used to raise indices as needed; in particular, with $g_{ij}^k = \sum_e
	g_{ije}\,g^{ek}$, we have the multiplication formula
\begin{equation}\label{gijk}
\basis_i \cp \basis_j = \smallsum{k}{} g_{ij}^k \; \basis_k.
\end{equation} 
\end{blanko}

\begin{blanko}{Intersection polynomials of $X$.}
	For a generic even element $\coh y \in H$, (identified with its coordinates
	with respect to $\basis_0,\ldots,\basis_r$ as in \ref{Gamma}), let
	$\phi(\coh{y}) \in \Q[[\fat{y}]]$ be the generating function for the
	integrals on $X$,
	\begin{eqnarray*}
		\phi(\coh y) \;\df\; \otrace{\exp(\coh y)} 
		& = & \sum_{n\geq 0} \frac{1}{n!}\otrace{\coh y^n} .  
	\end{eqnarray*}
	Note that the expansion is $\exp(y_0)$ times a polynomial
	in $y_1,\ldots,y_r$.

   Now set 
$$
   \phi_{ij} \df \frac{\d^2}{\d y_i \d y_j} \, \phi 
	= \otrace{\exp(\coh y)\cp\basis_i\cp\basis_j}
$$
and use the matrix $(g^{ef})$ to raise indices, putting
\begin{equation}\label{gandphi}
	 \phi^i{}_j \df \smallsum{e}{} g^{ie} \, \phi_{ej} ,
	   \quad\text{ }\quad 
	 \phi_j{}^i{} \df \smallsum{e}{}  \phi_{je}\, g^{ei} ,
	   \quad\text{ and }\quad 
	 \phi^{ij} =\sum_{e,f} g^{ie} \, \phi_{ef} \, g^{fj}.
\end{equation}
The two first entities appear in the following multiplication formula,
\begin{eqnarray}
	\exp(\coh y) \cp \basis_p \label{exptophi}
	&=& \smallsum{e,f}{} \, \basis_e \, g^{ef} \, \otrace{\basis_f \cp
\exp(\coh y) \cp \basis_p}  \\
	 & = & \smallsum{e}{} \, \basis_e \, \phi^e{}_p(\coh y),\notag
\end{eqnarray}
and similarly
$
\basis_q \cp \exp(\coh y) = \smallsum{f}{} \, \phi_q{}^f(\coh y) \, \basis_f .
$
(These are identities in $H[[\fat y]]$.)
\end{blanko}

\begin{lemma}
   \textsc{Sum formula.} Let $\coh y'$ and $\coh y''$ be generic even elements 
   of $H$. Then
   $$
    \phi(\coh y'+\coh y'') =
  \sum_{e,f} \phi_e(\coh y') \, g^{ef} \, \phi_f(\coh y'').
   $$
\end{lemma}
\begin{dem}
   Write down the exponential series for $X\times X = X'\times X''$,
   $$
   \exp\big( \smallsum{}{} y'_i \basis_i' + \smallsum{}{} y''_i \basis_i''\big)
   $$
	and integrate over the diagonal $\Delta$.  Computing this integral using the
	isomorphism $X\simeq \Delta$ yields the left hand side.  On the other hand,
	computing the integral using the K\"unneth formula for $\Delta$ and the
	projection formula yields the right hand side.
\end{dem}
This sum formula (and its derivatives) ensures that the expected
index-raising rules hold.  For example, the useful formula
\begin{equation}\label{phiijk}
   \phi^i{}_k{}^j = \smallsum{\ell}{} \, g_{k\ell}^i \, \phi^{\ell j}
\end{equation}
amounts to $\phi^i{}_k{}^j( \coh 0 + \coh y) =
	\smallsum{\ell,m}{}  \phi^i{}_{k\ell}(\coh 0)\,  g^{\ell m} \, 
	\phi_m{}^j(\coh y)
$.

\begin{blanko}{The deformed metric.}\label{defmet}
   Instead of using the integral $\otrace{} : H \pil \Q$, the new metric is based
   on the linear map $\ytrace{ } :H \pil \Q[[\fat{y}]]$ defined as
\begin{eqnarray*}
	\ytrace{\coh z} & \df & \otrace{ \exp(-2\coh y) \cp \coh z } .
\end{eqnarray*}
   It is thought of as a deformation of $\otrace{}$ since we recover this
   map upon setting $\coh y = \coh 0$.
   Now define the new metric $(\gamma_{ij})$ by
   $$
   \gamma_{ij} \df \gamma_{ij}(\coh y) \df
   \ytrace{\basis_i \cp \basis_j} = \phi_{ij}(-2\coh{y}),
   $$ 
	and adopt the obvious notation $\gamma_{ijk} = \ytrace{\basis_i \cp \basis_j\cp
	\basis_k} = \phi_{ijk}(-2\coh{y})$.  Let $(\gamma^{ij})$ denote the
	inverse matrix to $(\gamma_{ij})$.  Then it follows readily from the sum
	formula above that
   \begin{equation}  \label{gamma-def}
      \gamma^{ij} = \phi^{ij}(2\coh{y})  = \sum_{e,f}
      \phi^i{}_e \, g^{ef} \, \phi_f{}^j.
   \end{equation}
   The sum formula also yields
	$\sum_{e} \gamma_{ije} \, \gamma^{ef}
	= \sum_e \phi_{ije}(-2\coh{y})\; \phi^{ef}(2\coh{y}) =
   \phi_{ij}{}^f(\coh{0}) = g_{ij}^f$, so the intersection product can be
   written
   $$
   \basis_i \cp \basis_j = \sum_{e,f} \gamma_{ije} \, \gamma^{ef} \, \basis_f .
   $$
\end{blanko}

\begin{eks}\label{phiP2}
   For $\P^2$ (with $h\df c_1(\OO(1))$ and basis $\basis_i\df h^i$), we get
   $$
   (\gamma^{ij}) = \exp(2y_0)\begin{pmatrix}
   \phantom{0}0\phantom{0} & \phantom{0}0\phantom{0}& \phantom{0}1\phantom{0} \\
   0 & 1 & 2y_1 \\
   1 & 2y_1 & 2y_1^2 + 2y_2
   \end{pmatrix}.
   $$
   This matrix (with $y_0$ and $y_2$ set to zero) was first written down in the 
   pioneer
   article~\cite{diFrancesco-Itzykson} of Di~Francesco and Itzykson.
\end{eks}

\section{WDVV equations and the tangency quantum product}

\begin{blanko}{Notation.}
	Let lower indices denote partial derivatives, e.g.
	$$
	\Gamma_{ij} \df \Gamma_{x_i x_j} \df \frac{\d^2}{\d x_i \d x_j}\Gamma =
	\sum_{\beta>0} q^\beta 
	\brac{ \exp( \mtau_0(\coh x) + \mtau_1(\coh y))\cdot\mtau_0(\basis_i)\,
	\mtau_0(\basis_j)}\!_\beta ,
	$$
	and set
\begin{equation}  \label{Gammaij}
   \Gamma_{(ij)} = \Gamma_{(x_i x_j)} \df \sum_{k=0}^r g_{ij}^k \ \Gamma_{x_k},
\end{equation}
--- the ``directional derivative with respect to the product
$\basis_i\cp \basis_j = \sum g_{ij}^k \basis_k$''.
\end{blanko}

\begin{satz}\label{WDVV}
	The following form of the WDVV equations holds for the tangency quan\-tum
	potential.
\begin{eqnarray*}
	\Gamma_{(ij)k\ell} + \Gamma_{ij(k\ell)} + 
	\sum_{e,f}\; \Gamma_{ije} \ \gamma^{ef} \ \Gamma_{fk\ell} & = & \pm\big(
   \Gamma_{(jk)i\ell} + \Gamma_{jk(i\ell)} + 
	\sum_{e,f}\; \Gamma_{jke} \ \gamma^{ef} \ \Gamma_{fi\ell}\big),
\end{eqnarray*}
where $\pm$ denotes the sign $(-1)^{\deg \basis_i(\deg \basis_j +\deg\basis_k)}$.
\end{satz}
\begin{obs}
   In the special case $X=\P^2$, it turns out the only non-trivial relation is
   the one with $i=j=1$ and $k=\ell=2$.  Since $\Gamma_0 = 0$ by the string 
	equation (cf.~\cite{Grab-Kock-Pand}), the relation then reads
$$
\Gamma_{222} = \exp(2y_0)\big(\Gamma_{112}^2 - \Gamma_{111}\Gamma_{122} + 
2y_1 ( \Gamma_{122}\Gamma_{112} - \Gamma_{111}\Gamma_{222} )
+ (2y_1^2 + 2y_2)( \Gamma_{122}^2 - \Gamma_{112}\Gamma_{222})\big).
$$
This equation was first found by L.~Ernstr\"om and G.~Kennedy~\cite{EK2},
cf.~Remark~\ref{EKremark} below, while the special case of $y_0=y_2=0$ goes back
to Di~Francesco-Itzykson~\cite{diFrancesco-Itzykson}.  Setting also $y_1=0$ we
are back to the celebrated formula of M.~Kontsevich~\cite{KM:9402}.
\end{obs}
\begin{dem}
	The proof follows the line of arguments of the proof of the WDVV equations
	for the usual Gromov-Witten potential, cf.~Kontsevich-Manin~\cite{KM:9402}. 
	The novelty is the splitting lemma for enumerative descendants and the
	appearance of the deformed metric.
	
	For fixed $\beta>0$ and $a,b\geq 0$, consider the moduli stack $\ov
	M_{0,\{p_1,p_2,p_3,p_4\}\cup S}(X,\beta)$ where $S$ is a marking set of
	cardinality $a+b$.  Consider the product
	$$
	\mtau_0(\basis_i)\mtau_0(\basis_j)\mtau_0(\basis_k)\mtau_0(\basis_\ell) \ 
	\frac{\mtau_0(\coh x)^a}{a!}\frac{\mtau_1(\coh y)^b}{b!}
	$$
	where the first four classes correspond to the marks $p_1,p_2,p_3,p_4$.
	Now integrate the product over each side of the
	fundamental equivalence
   \begin{equation}\label{equiv}
   (p_1, p_2 \mid p_3, p_4) = ( p_2, p_3 \mid p_1, p_4),
   \end{equation}
	where $(p_1, p_2 \mid p_3, p_4)$ denotes the sum of all boundary divisors
	having $p_1$ and $p_2$ on one twig and $p_3$ and $p_4$ on the other.  Summing
	up these equations over all $a,b\geq 0$ and over all $\beta>0$ (as in the
	definition of $\Gamma$) we'll get the desired equation.  Let us treat the
	left hand side.  On the right hand side of the equation the arguments are the
	same; only it is necessary initially to permute the four special factors,
	which accounts for the sign.

	On the left hand side of the equation, let us first consider the contribution
	from the trivial degree partitions, say $\beta'=0$.  Then the only possible
	distribution of the marks giving contribution is when all the spare marks
	fall on the right hand twig, which leaves us with the single boundary divisor
	$D_{12}$.  Now according to \ref{D_ij}, the effect of multiplication with
	this divisor is to merge the two classes $\mtau_0(\basis_i)$ and
	$\mtau_0(\basis_j)$, so in the end we get
   $$
   \brac{\mtau_0(\basis_i\cp \basis_j)\mtau_0(\basis_k)\mtau_0(\basis_\ell)
\frac{\mtau_0(\coh x)^a}{a!}\frac{\mtau_1(\coh y)^b}{b!}}\!_{\beta} .
$$
	Similarly, the case $\beta''=0$ gives
	$\brac{\mtau_0(\basis_i)\mtau_0(\basis_j)\mtau_0(\basis_k\cp \basis_\ell)
	\frac{\mtau_0(\coh x)^a}{a!}\frac{\mtau_1(\coh y)^b}{b!}}\!_{\beta}$. 
	Summing over $\beta>0$ and $a,b\geq 0$ gives exactly the two linear terms on
	the left hand side of the promised equation, cf.~(\ref{Gammaij}).
      
   Now for those boundary divisors in the linear equivalence corresponding to
   strictly positive degree partitions.  To each irreducible component on the
   left hand side of (\ref{equiv}), we apply the Splitting Lemma~\ref{split2},
   getting all together
	   \begin{multline*}
   \bigsum{}{}\sum_{p,q}
   \brac{\frac{\mtau_{0}(\coh x)^{a'}}{a'!}\frac{\mtau_1(\coh 
   y)^{b'}}{b'!}\mtau_0(\basis_i)\mtau_0(\basis_j) \,
   \mtau_0(\frac{\coh y^{c'}}{c'!}\!\cp \basis_p)}\!_{\beta'}
   \\ g^{pq} \  
   \brac{\mtau_0(\basis_q\cp\frac{\coh y^{c''}}{c''!})
   \,\mtau_0(\basis_k)\mtau_0(\basis_{\ell})\frac{\mtau_{0}(\coh x)^{a''}}{a''!}
	\frac{\mtau_1(\coh y)^{b''}}{b''!}}\!_{\beta''}.
   \end{multline*}
where the big outer sum is over all $a'+a''=a$ and all $b'+c' + b''+c''= b$.
(The $b'$ corresponds to what was called $b'-s'$ in the splitting lemma.)
Now sum over all $a$ and $b$ getting
   \begin{multline*}
   \sum_{p,q} 
   \brac{\exp(\mtau_{0}(\coh x)+\mtau_1(\coh 
   y))\cdot\mtau_0(\basis_i)\mtau_0(\basis_j) \,
   \mtau_0(\exp(\coh y)\!\cp \basis_p)}\!_{\beta'}
   \\ g^{pq} \  
   \brac{\mtau_0(\basis_q\cp\exp(\coh y))
   \,\mtau_0(\basis_k)\mtau_0(\basis_{\ell})\cdot\exp(\mtau_{0}(\coh x)+
	\mtau_1(\coh y))}\!_{\beta''}.
   \end{multline*}
	Next, use equation (\ref{exptophi}) to get rid of $\exp(\coh y)$,
	and sum over all $\beta>0$ as in the definition of $\Gamma$, arriving at
$$
\underset{e,f}{\sum_{p,q}} \ \Gamma_{ije} \, \phi^e{}_p \, g^{pq} \, 
\phi_q{}^f \, \Gamma_{fk\ell} .
$$
By equation~(\ref{gamma-def}),  this is just the quadratic term of the left hand
side of the desired equation.
\end{dem}

\begin{blanko}{Topological recursion relation.}
	The above deformed WDVV equations alone are not sufficient to determine all
	the first enumerative descendants from the primary Gromov-Witten invariants.
	But topological recursion is available also for the enumerative descendants
	cf.~\cite{Grab-Kock-Pand}.  In the present set-up, that topological recursion
	relation takes the following pleasant form.
   $$
   \Gamma_{y_i x_j x_k } = \Gamma_{x_i(x_j x_k)} - \Gamma_{(x_i x_j) x_k}
   - \Gamma_{(x_i x_k) x_j} +
   \sum_{e,f} \; \Gamma_{x_i x_e} \ \gamma^{ef} \ \Gamma_{x_f x_j x_k}.
   $$
The shape of this equation stems from the boundary expression of the modified psi
class, $\mpsiclass_i = (p_i\mid p_j, p_k) - \boundary_i$, where $\boundary_i$ is
the sum of all boundary divisors such that $p_i$ is on a contracting twig.  
The deformed metric enters in the quadratic terms for the same reason as in the 
proof above.
\end{blanko}

\begin{blanko}{The tangency quantum product.}\label{tprod}
	The {\em tangency quantum product} is the $\ring[[\fat{x},\fat{y}]]$-bilinear
	product on $H\upperstar (X,\Q) \tensor_\Q \ring[[\fat{x},\fat{y}]]$ defined
	by the rule
$$
\basis_i * \basis_j \df \basis_i \cp \basis_j + 
\sum_{e,f} \; \Gamma_{ije} \ \gamma^{ef} \ \basis_f.
$$
\end{blanko}
Observe that this product specialises to the usual quantum product upon setting 
the formal variables $\fat{y}$ to zero.  Clearly the product is supercommutative.

\begin{satz}\label{ass}
   The tangency quantum product is associative.
\end{satz}
\begin{dem}
   This is a straightforward consequence of Theorem~\ref{WDVV};  it amounts to 
   checking the associativity relations on the generators, using the definition
   of the product.  The only subtle point in the verification is the identity
   $$
      \sum_{e,f} \;\Gamma_{ije} \ \gamma^{ef} \ (\basis_f \cp \basis_k) 
   = \sum_{\ell,m}\;\Gamma_{ij(k\ell)} \ \gamma^{\ell m} \ \basis_m ,
$$
which follows from the properties of the structure constants $g_{ij}^k$
listed in (\ref{gijk}), (\ref{phiijk}), and (\ref{Gammaij}).
\end{dem}

\begin{obs}\label{EKremark}
	In the special case $X=\P^2$, this product (or a simple change of coordinates
	of it) was previously constructed via ad hoc methods by L.~Ernstr\"om and
	G.~Kennedy~\cite{EK2}, who also gave a tour de force proof of its
	associativity.  Their construction relies on the space of stable lifts, and
	seems to be peculiar to the projective plane.
\end{obs}

\bigskip

\section{Integrability}
\label{integrability}

\begin{blanko}{The classical potential} for $X$ is the generating function for
the triple top products
\begin{eqnarray*}
\otrace{ \frac{\coh x^3}{3!}}
&=&\sum_{i,j,k} \frac{x_k x_j x_i}{6} 
\otrace{ \basis_i \cp \basis_j \cp \basis_k }.
\end{eqnarray*}

By construction, its third derivatives are just $g_{ijk}$, the structure
constants for the cup multiplication (in the Poincar\'e metric).  In usual
quantum cohomology, this potential is reinterpreted as the $\beta=0$ part of the
Gromov-Witten potential: a quantum correction (the $\beta>0$ part) is added to
the classical potential, in such a way that the third derivatives of this sum
are the structure constants of a new associative product, the quantum product.
\end{blanko}

\begin{blanko}{``The tangency classical potential''.}\label{tcp}
The tangency quantum product relates to $\ytrace{}$ exactly as the usual 
quantum product relates to $\otrace{}$. So let us introduce a potential
\begin{eqnarray*}
\Phi^0(\coh{x},\coh{y}) \df \ytrace{\frac{\coh x^3}{3!}} 
= \sum_{i,j,k} \frac{x_k x_j x_i}{6} \ytrace{\basis_i
\cp \basis_j \cp \basis_k} ,
\end{eqnarray*}
which to fit into the picture could be called the tangency classical potential,
although it is neither classical, nor has anything particular to do with
tangency.  By construction, its third derivatives are $\gamma_{ijk}$, the
structure constants of the cup multiplication, but this time in the deformed
metric:
$$
\basis_i \cp \basis_j = \sum_{e,f} \Phi^0_{ije} \ \gamma^{ef} \ \basis_f .
$$
\end{blanko}

\begin{blanko}{The tangency potential, including $\beta=0$.}\label{beta0}
   Introducing the potential
$$
\Phi(\coh{x},\coh{y}) \df \Phi^{0}(\coh{x},\coh{y}) + 
\Gamma(\coh{x},\coh{y}),
$$
whose third derivatives are $\Phi_{ijk} = \gamma_{ijk} + \Gamma_{ijk}$,
the tangency quantum product can be written
$$
\basis_i * \basis_j =
\sum_{e,f} \Phi_{ije} \ \gamma^{ef} \ \basis_f,
$$
and the  WDVV equation~\ref{WDVV} then takes the usual form
   $$
   \sum_{e,f}\; \Phi_{ije} \ 
   \gamma^{ef} \ \Phi_{fk\ell} 
   = 
   \pm \sum_{e,f}\; \Phi_{jke} \ 
   \gamma^{ef} \ \Phi_{fi\ell}.
   $$
\end{blanko}

\begin{blanko}{Interpretation of the degree-zero term.}
   \label{interpret0}
   Writing down this potential $\Phi=\Phi^0 + \Gamma$ calls for an
   interpretation of $\Phi^0$ in terms of some top products on the degree-zero
   spaces.  Unfortunately there is no way to define the modified psi class on
   $\ov M_{0,n}$, if we want it to satisfy the two rules: (i) it should be
   compatible with pull-back along forgetful morphisms, and (ii) it should
   satisfy the push-down formula $\pi\lowerstar \mpsiclass = -2$, independent of
   the number of marks, to give the dilaton equation for modified psi classes
   (cf.~\cite{Grab-Kock-Pand}). Basically this is impossible because the
   1-pointed spaces do not exist in degree zero.
   The best one can do is to define the class on a fixed 4-pointed space and
   then pull it back to the hierarchy lying over this space, but this definition
   depends on the choice of the three extra marks.  
   
So take a moduli space $\ov M_{0,S\cup \{q_1,q_2,q_3\}}(X,0)$ with three
distinguished marks.  For each of the other marks $p_s\in S$, define the
modified psi class $\mpsiclass_s$ as the pull-back from $\ov
M_{0,\{p_s,q_1,q_2,q_3\}}\simeq \P^1$ of the class of degree $-2$.  One easily
checks that this is equivalent to defining $\mpsiclass_s \df \psiclass_s -
\boundary_s$, where $\boundary_s$ is the sum of all boundary divisors such that
$p_s$ is on a (contracting) twig together with at most one of the distinguished
marks.  (This description is then
compatible with the boundary description of the modified psi class in the
$\beta>0$ case (cf.~\cite{Grab-Kock-Pand}), since in that case there are no
distinguished marks.)

In this setting, define the invariant 
$$
\brac{\frac{\mtau_0(\coh x)^a}{a!}\frac{\mtau_1(\coh y)^b}{b!}\,
\mtau_0(\basis_i)\mtau_0(\basis_j)\mtau_0(\basis_k)}\!_0
$$
in the obvious way, with the last three classes corresponding to the three
distinguished marks.  Identifying $\ov M_{0,S\cup\{q_1,q_2,q_3\}}(X,0)$ with
$\ov M_{0,S\cup\{q_1,q_2,q_3\}} \times X$, all the evaluation morphisms are just
the projection $p$ to $X$, so the integrand has a factor $p\upperstar(\frac{\coh
x^a}{a!}\cp\frac{\coh y^b}{b!} \cp\basis_i\cp\basis_j\cp\basis_k)$.  The
remaining factors are modified psi classes from $\ov
M_{0,S\cup\{q_1,q_2,q_3\}}$; since each of them is alone on its mark, we can
push them down one by one, arriving at a factor $(-2)^b$.  The remaining
integral $\int 1$ is zero for dimension reasons, unless we have come down to
just $\ov M_{0,\{q_1,q_2,q_3\}}$, which means $a=0$.  We conclude
$$
\brac{\frac{\mtau_1(\coh x)^b}{b!}\,
\mtau_0(\basis_i)\mtau_0(\basis_j)\mtau_0(\basis_k)}\!_0
=
\otrace{ \frac{(-2\coh y)^b}{b!} \cp \basis_i\cp \basis_j
\cp \basis_k }.
$$
Summing over all $b$ (and $a$) we get
$$
\otrace{ \exp(-2\coh y) \cp \basis_i \cp \basis_j \cp \basis_k} 
= \ytrace{\basis_i \cp \basis_j \cp \basis_k}
= \gamma_{ijk}
$$
showing that at least the third derivatives of $\Phi^0$ have an 
interpretation as top products on degree-0 moduli spaces, and after all
it is the third derivatives that really matter.
\end{blanko}

Once we know that the structure constants of the tangency quantum product are 
third derivatives of the single potential $\Phi$, we are in position to 
given an interpretation in terms of Frobenius manifolds.

\begin{blanko}{Formal Frobenius manifolds.}
	For convenience let us recall (from Manin~\cite{Manin}, Ch.~III) the
	definition of a formal Frobenius manifold (over $k$).  Let $k$ be a
	supercommutative $\Q$-algebra.  Let $H$ be a free $k$-module of finite rank,
	with generators $\basis_0,\ldots,\basis_r$, and let $g : H \tensor H \pil k$
	denote an even symmetric non-degenerate bilinear pairing.  Let $K=k[[H^t]]$
	be the completed symmetric algebra of the dual module $H^t$.  In other words,
	if
%
	$\coh x = \sum x_i \basis_i$ is a generic even element of $H$ 
	then $K=k[[x_0,\ldots,x_r]]$.  Now the structure of a formal Frobenius
	manifold on $(H,g)$ is given by an even potential $\Phi\in K$ (defined up to
	quadratic terms) satisfying WDVV. In other words, the multiplication
	$\basis_i * \basis_j \df \sum_{e,f} \Phi_{ije}\; g^{ef}\; \basis_f$ makes
	$H\tensor_k K$ into an associative supercommutative $K$-algebra.
\end{blanko}

\smallskip

With $H\df H\upperstar (X,\ring)$, the results of this section readily imply:

\begin{satz}
	The cohomology $\ring[[\fat{y}]]$-module $H[[\fat{y}]]$ with bilinear
	non-degenerate pairing $\gamma : H[[\fat{y}]] \tensor H[[\fat{y}]] \pil
	\ring[[\fat{y}]]$, equipped with the tangency quantum potential $\Phi \in
	\ring[[\fat{x},\fat{y}]]$ constitutes a formal Frobenius manifold over
	$\ring[[\fat{y}]]$.  \qed
\end{satz}

In fact, this formal Frobenius manifold is a deformation over $\Q[[\fat{y}]]$
of the formal Frobenius manifold of usual quantum cohomology.
While the underlying space is trivially deformed under
this deformation, the metric and the potential vary non-trivially.


\begin{thebibliography}{10}

\bibitem{Behrend}
{\sc Kai Behrend}.
\newblock {\em {G}romov-{W}itten invariants in algebraic geometry}.
\newblock Invent. Math. {\bf 127} (1997), 601--617.
\newblock (alg-geom/9601011).

\bibitem{Behrend-Manin}
{\sc Kai Behrend {\rm and }Yuri~I. Manin}.
\newblock {\em Stacks of stable maps and {G}romov-{W}itten invariants}.
\newblock Duke. J.~Math. {\bf 85} (1996), 1--60.
\newblock (alg-geom/9506023).

\bibitem{diFrancesco-Itzykson}
{\sc Philippe di~Francesco {\rm and }Claude Itzykson}.
\newblock {\em Quantum intersection rings}.
\newblock In {R. Dijkgraaf, C. Faber,} {\rm and }G.~{van der Geer}, editors,
  {\em The moduli space of curves}, vol. 129 of Progress in Mathematics, 
  pp.~81--148. Birkh{\"a}user, Boston, MA, 1995.

\bibitem{EK2}
{\sc Lars Ernstr{\"o}m {\rm and }Gary Kennedy}.
\newblock {\em Contact cohomology of the projective plane}.
\newblock Amer. J.~Math. {\bf 121} (1999), 73--96.
\newblock (alg-geom/9703013).

\bibitem{FP-notes}
{\sc William Fulton {\rm and }Rahul Pandharipande}.
\newblock {\em Notes on Stable Maps and Quantum Cohomology}.
\newblock In {J.~Koll{\'a}r, R. Lazarsfeld} {\rm and }D.~Morrison, 
editors,\linebreak
  {\em Algebraic Geometry, Santa Cruz 1995}, vol. 62, II of Proc. Symp. Pure.
  Math., \linebreak pp.~45--96.
\newblock (alg-geom/9608011).

\bibitem{Getzler:9612}
{\sc Ezra Getzler}.
\newblock {\em Intersection theory on {$\overline M_{1,4}$} and elliptic
  {G}romov-{W}itten invariants}.
\newblock J.~Amer. Math. Soc. {\bf 10} (1997), 973--998.
\newblock (alg-geom/9612004).

\bibitem{Grab-Kock-Pand}
{\sc Tom graber, Joachim Kock, {\rm and }Rahul Pandharipande}.
\newblock Descendant invariants and characteristic numbers (in genus $0$, $1$, 
and $2$).
\newblock In preparation.

\bibitem{_9902}
{\sc Joachim Kock}.
\newblock Recursion for twisted descendants and characteristic numbers of
  rational curves.
\newblock Preprint, math.AG/9902021.

\bibitem{_thesis}
{\sc Joachim Kock}.
\newblock {\em Tangency quantum cohomology and enumerative geometry of 
\linebreak 
rational curves}.
\newblock PhD thesis, Recife, Brazil, March 2000.
\newblock Available at \\ http://www.math.kth.se/{\~{}}kock/tese/tese.ps.

\bibitem{KM:9402}
{\sc Maxim Kontsevich {\rm and }Yuri~I. Manin}.
\newblock {\em {G}romov-{W}itten classes, quantum cohomology, and enumerative
  geometry}.
\newblock Comm. Math. Phys. {\bf 164} (1994), 525--562.
\newblock (hep-th/9402147).

\bibitem{Li-Tian}
{\sc Jun Li {\rm and }Gang Tian}.
\newblock {\em Virtual moduli cycles and {G}romov-{W}itten invariants of
  algebraic varieties}.
\newblock J.~Amer. Math. Soc. {\bf 11} (1998), 119--174.
\newblock (alg-geom/9602007).

\bibitem{Manin}
{\sc Yuri~I. Manin}.
\newblock {\em {F}robenius manifolds, quantum cohomology, and moduli spaces}.
\newblock AMS Colloquium Publications, Providence, RI, 1999.

\bibitem{Pand:mail}
{\sc Rahul Pandharipande}.
\newblock Unpublished e-mail to {L}ars {E}rnstr{\"o}m, {M}ay 2, 1997.

\bibitem{Witten}
{\sc Edward Witten}.
\newblock {\em Two-dimensional gravity and intersection theory on moduli
  space}.
\newblock Surveys in Diff. Geom. {\bf 1} (1991), 243--310.

\end{thebibliography}

\end{document}